\documentclass[12pt]{article}
\usepackage{amssymb}
\usepackage{amsthm}\usepackage{epsf,epsfig,amsfonts,}
\usepackage{amsmath} \usepackage[french]{babel}
\usepackage{graphicx} 

\usepackage{a4wide}

\newtheorem{corollary*}{Corollary}

\newcommand{\be}{\begin{equation}}
\newcommand{\ee}{\end{equation}}

\newcommand{\weg}[1]{}

\usepackage{epsf,epsfig,amsfonts,a4wide}

\usepackage{amsfonts}
\usepackage{amsmath,amssymb,amsthm,color}
\usepackage{url}

\newtheorem{Th}{Theorem}

\theoremstyle{remark}

\newtheorem{Rem}{Remark}

\newcommand{\const}{\mbox{\rm const}}
 \title{There are no conformal Einstein rescalings of complete pseudo-Riemannian Einstein metrics
}
\date{} \author{ Volodymyr  Kiosak, Vladimir  S. Matveev\thanks{ Institute of Mathematics, FSU Jena, 07737 Jena Germany,  vladimir.matveev@uni-jena.de}}

\begin{document}
\maketitle 

\begin{abstract}
{\bf Il n'existe pas de rescaling conform\'ement Einstein
       d'une m\'etrique d'Einstein pseudo-Riemannienne compl\`ete.}

Soit $g$ une m\'etrique pseudo-riemannienne non d\'efinie de type Einstein  telle que la
m\'etrique conform\'ement \'equivalente $\psi^{-2} g$ soit aussi d'Einstein.
Nous montrons que si  la m\'etrique $g$ est lumi\`ere-compl\`ete, i.e. ses g\'eod\'esiques isotropes sont compl\`etes,    alors le coefficient $\psi$ est constant.  Si la vari\'et\'e est
ferm\'ee, l'hypoth\`ese de  compl\'etude peut \^etre omise (ce dernier
r\'esultat est d\^u \`a Mikes-Radulovich et K\"uhnel, mais notre
d\'emonstration est plus courte).

La d\'emonstration est bas\'ee sur l'\'etude du comportement de la
fonction $\psi$ le long des g\'eod\'esiques de  type lumi\`ere. Si   
   $\gamma(t)$ est une telle g\'eod\'esique, alors: 
  $\psi(\gamma(t))= \const_1 \cdot  t  + \const_2$.  Comme  la fonction
$\psi$ est non-nulle,  la lumi\`ere-compl\'etude   implique $\psi =\const_2$.

Si la vari\'et\'e est ferm\'ee, la fonction $\psi$ prend sa valeur maximale
$\psi_{max}$ en un certain  point.  Donc, pour toute g\'eod\'esique
de type  lumi\`ere $\gamma$ passant par ce point, on a  $\const_1=0$,
ce qui implique que  $\psi=\psi_{max}$ en tout  point de cette g\'eod\'esique.
En r\'ep\'etant cet argument, on obtient  que pour toute
g\'eod\'esique $\gamma_1$  de type lumi\`ere coupant $\gamma$, 
$\psi=\psi_{max}$ en tout  point de $\gamma_1$, et    ainsi de suite.
On en d\'eduit $\psi$ est constante sur la
vari\'et\'e enti\`ere, car   deux points quelconques  peuvent   \^etre joints  par une suite de
g\'eod\'esiques   de type lumi\`ere,.
 \\[.1cm]

\centerline{\bf Abstract}

\vspace{1ex} 
Let $g$ be an Einstein metric of indefinite signature such that the conformally-equivalent metric 
$\psi^{-2} g$ is also  Einstein. We show that if  the metric $g$ is light-line complete, then the  conformal coefficient $\psi$ is constant. 
If the manifold is closed, the completeness assumption   can be omitted (the latter  result is due to Mikes-Radulovich and K\"uhnel, but our proof is much simpler).

The proof is  based on the investigation of the   behavior of the function $\psi$ along light-line geodesics: we show that for every light-line geodesic $\gamma(t)$ we have $\psi(\gamma(t))= \const_1 \cdot  t  + \const_2$. Since the function $\psi$ can not vanish, the light-line completeness of the metric implies  $\psi =\const_2$. 

 If the manifold is closed, the function $\psi$  accepts its maximal value $\psi_{max}$  at a certain point. Then, for every light-line geodesic $\gamma$ through this point we have $\const_1=0$ implying $\psi=\psi_{max}$ at every point of this geodesic. Repeating the argumentation, we obtain that for every  light-line geodesic $\gamma_1$ intersecting $\gamma$ we have $\psi=\psi_{max}$ at every point of $\gamma_1$ as well and so on. Since every two points can be connected by a sequence of light-line geodesics,  $\psi$  is constant on the whole manifold.

\end{abstract} 

\begin{Th}  \label{einstein1}
Let $g$ be a light-line-complete pseudo-Riemannian Einstein metric of indefinite signature (i.e., for no constant $c$ the metric $c\cdot g$ is Riemannian) on a connected $(n>2)-$dimensional manifold $M$.    Assume that for the nowhere vanishing function $\psi$  the metric    $\psi^{-2} g$ is also  Einstein. Then,   $\psi$ is a constant.
\end{Th}

\begin{Rem} {Theorem}  \ref{einstein1}   fails  for  Riemannian metrics  (even if we replace light-line completeness by usual  completeness) -- M\"obius transformations of the standard round sphere and the stereographic map of the punctured  sphere to the  Euclidean space are conformal nonhomothetic mappings.    One can  construct other  
examples on warped Riemannian manifolds, see  \cite[Theorem 21]{kuehnel}. \end{Rem} 
\begin{Rem} By Theorem \ref{einstein1},  {\it  light-line  complete pseudo-Riemannian Einstein metrics  of indefinite signature do not admit nonhomothetic  conformal  complete vector fields}.
  The Riemannian version of this result is due to Yano and Nagano   \cite{nagano}. Moreover, the assumption that the metric is Einstein can be omitted (by the price of considering only essential conformal vector fields):  as it was proved by  D. Alekseevksii~\cite{Al},    J. Ferrand~\cite{Fe2} and R. Schoen~\cite{schoen},   a  Riemannian manifold   admitting an    essential   complete vector field  is  conformally equivalent to  the round sphere or  the Euclidean space.  It is still not known whether the last statement  (sometimes called Lichnerowicz-Obata conjecture) can  be extended to the pseudo-Riemannian case, see \cite{leitner} for a  counterexample in the $C^1 -$ smooth category, and \cite{frances,frances2} for a good survey on this topic.  \end{Rem} 
\begin{Rem} In the 4-dimensional lorenz case, Theorem \ref{einstein1} was  known in folklore:   more precisely, 
conformal Einstein rescalings of 4-dimensional  Einstein metrics were   described by Brinkmann \cite{brinkmann}, see also \cite[Corollary 2.10]{kuehnel2}. The list of all such metrics and their conformal Einstein rescalings is pretty simple   and one can 
directly verify our Theorem \ref{einstein1} by calculations.
\end{Rem} 

\begin{Rem}  A partial case of Theorem \ref{einstein1} is   \cite[Theorem 2.2]{kuehnel2}, in which it is  assumed that   both metrics are complete. This extra-assumption is very natural in the context of  \cite{kuehnel2} since the paper is dedicated to the classification of conformal vector fields; moreover, Theorem 2.2 is not the main result of the paper. It  is not clear whether in the proof of  \cite[Theorem 2.2]{kuehnel2}  the assumption that the second  metric is  complete could be omitted. \end{Rem}

{\bf Proof of Theorem \ref{einstein1}. } It is well-known (see for example \cite[eq. (2.21)]{brinkmann} or \cite[Lemma 1]{kuehnel}) that the  Ricci curvatures $R_{ij}$ and $\bar R_{ij}$ of two conformally equivalent metrics $g$ and $\bar g= \psi^{-2} g = e^{-2\phi} g$  are related  by 
\begin{equation} \label{eq1}  
\bar R_{ij} =R_{ij}+ (\Delta \phi- (n-2) \|\nabla \phi\|^2)g_{ij}+ \frac{n-2}{\psi}\nabla_i\nabla_j \psi.   
\end{equation} 

Consider a light-line geodesic $\gamma(t)$ of the metric $g$. Since the metric $g$ is light-line-complete, $\gamma(t)$ is defined on the whole $\mathbb{R}$. ``Light-line" means that  $g(\dot\gamma(t), \dot\gamma(t)) = g_{ij}\dot \gamma^i(t) \dot \gamma^j(t)=0$, where $\dot\gamma$ is the velocity vector of $\gamma$ (it is well-known that if this property is fulfilled in one point then it is fulfilled at every point of the geodesic).  

Now contract  \eqref{eq1} with $\dot \gamma^i \dot \gamma^j$.
 Since the metrics are Einstein and conformally equivalent, 
  $\bar R_{ij}$, $R_{ij}$ and  $g_{ij}$ are  proportional to $ g_{ij}$,  and therefore  the only term which does not vanish is $\dot \gamma^i \dot \gamma^j\frac{n-2}{\psi}\nabla_i\nabla_j \psi$.  Thus, $\dot \gamma^i \dot \gamma^j\nabla_i\nabla_j \psi=0$. 
 
 Clearly, at every point of the geodesic we have $\dot \gamma^i \dot \gamma^j\nabla_i\nabla_j \psi= \tfrac{d^2}{dt^2}\psi(\gamma(t))$. Thus, $ \tfrac{d^2}{dt^2}\psi(\gamma(t))=0$ implying $ \psi(\gamma(t))= \const_1 \cdot t + \const$. Since by assumptions the function $\psi$ is defined on the whole $\mathbb{R}$ and  is equal  to zero at no point, we have $\const_1=0$ implying $\psi\equiv \const$ along every light-line geodesic. 

Now, 
 every two points of a connected manifold can be connected by a finite  sequence  of light-line geodesics.
 Indeed,  consider $\mathbb{R}^n $  with the standard pseudo-Euclidean  metric $g_0$  of  the same 
 signature $(r, n-r)$, $1\le r  < n$ as the metric $g$.  The union  of  all light-line geodesics passing  through points $a $  (resp. $b$) are the 
   standard  cones   $C_a:= \{(x_1,...,x_n)\in \mathbb{R}^n \mid (x_1- a_1)^2+...+(x_r-a_r)^2-(x_{r+1}- a_{r+1})^2-...-(x_n-a_n)^2=0\}$  and, resp.,   $C_b:= \{(x_1,...,x_n)\in \mathbb{R}^n \mid (x_1- b_1)^2+...+(x_r-b_r)^2-(x_{r+1}- b_{r+1})^2-...-(x_n-b_n)^2=0\}$. These two cones    always have points of transversal  intersection.  Thus, two arbitrary points of $\mathbb{R}^n$ can be connected by a sequence of two light-line geodesics of $g_0$. Since the restriction of the metric $g$ to  a small neighborhood $U\subseteq M^n$  can be viewed as a small perturbation of the  metric $g_0$   in $\mathbb{R}^n$, two points in $U$ can be connected by a sequence of two light-line geodesics.  Then, the set of points of $M$ that can be connected with  a fixed  point $p\in M^n$ by a finite sequence of light-line geodesics  is open and  closed implying it coincides with $M$.

 Since every two points of $M$  can be connected by a sequence  of light-line geodesics, and since  as we proved above the function $\psi$  is constant along every light-line geodesic, we have that $\psi$ is constant on the whole manifold as we claimed, \qed

\begin{Th}
\label{einstein2}
Let $g$ be a  pseudo-Riemannian Einstein metric of indefinite signature  on a connected closed (i.e., compact with  no boundary)  $(n>2)-$dimensional manifold $M$.    Assume that for the nowhere vanishing function $\psi$  the metric    $\psi^{-2} g$ is also  Einstein. Then,   $\psi$ is a constant.
\end{Th} 

\begin{Rem}
Theorem \ref{einstein2} is not new and is   in  \cite[Theorem 5]{mi}. Moreover, Wolfgang K\"uhnel explained us  how one can obtain the proof combining the results of PhD thesis  of  Kerckhove,  equation \eqref{eq1} due to \cite{brinkmann},  and  also \cite[Proposition 3.8(1)]{kuehnel2}. 
Our  proof of   Theorem  \ref{einstein2} is much easier than the  proofs of Mikes--Radulovich and K\"uhnel.
Actually, the initial version  of our  paper  did not contain Theorem 2 at all, but  after J. Mikes sent us  his paper we immediately saw that the  proof  of their Theorem 5 can be essentially simplified by using the trick from the proof of our Theorem \ref{einstein1}. 
 \end{Rem}

 {\bf Proof of Theorem \ref{einstein2}. }  Since $M$ is closed, there exists $p_0\in M$ such that the value of $\psi$ is maximal (we denote this value by $\psi_{max}$). 
 We take a light-line geodesic  $\gamma$ such that $\gamma(0)= p_0$. As we explained in the proof of Theorem \ref{einstein1}, the function $\psi(\gamma(t))$ is equal to
 $\const \cdot t + \psi_{max}$. Since the value of $\psi$ at the point $p_0$ is maximal, $\const =0$ implying $\psi(\gamma(t))\equiv \psi_{max}$. Then, for every  point $p_1$  of geodesic $\gamma$ the value of $\psi$ is maximal.   We can therefore repeat the argumentation and show that   for every light-line geodesic $\gamma_1$ such that $\gamma_1(0)=p_1$ we have  $\psi(\gamma_1(t))\equiv \psi_{max}$  and so on.  Since every two points of $M$  can be connected by a sequence  of light-line geodesics, we have that $\psi$ is constant on the whole manifold, \qed  
 
{\em Acknowledgement:} We thank W. K\"uhnel and H.-B. Rademacher for sending us the preliminary version of their  paper \cite{kuehnel2} and for useful discussions,  and  J. Mikes for sending  us his  paper \cite{mi}.

   When we obtained the proof, we asked all experts we know whether the proof is new, and are  grateful to those who answered, in particular to M. Eastwood, Ch. Frances, R. Gover, G. Hall, F. Leitner,  and  P. Nurowski. 

We thank Konrad  Sch\"obel and Ghani Zeghib for helping us  with the translation of the title and the  abstract to French.   
    
     Both authors were  partially supported by Deutsche Forschungsgemeinschaft (Priority Program 1154 --- Global Differential Geometry), and by  FSU Jena.

\end{document}